%% file: biflat.tex
\documentclass[11pt]{article}
\input{format}
\input{mathdefs}

\input{thmdefs}

\title{Operator biflatness of the Fourier algebra \\
and approximate indicators for subgroups}
\author{{\it Oleg Yu.\ Aristov}\thanks{Work on this paper was done when the author visited the universities of Alberta and of Waterloo in the fall of 2001. Financial support through the NSERC grants no.\ 227043-00 and
no.\ 90749-00 is gratefully acknowledged.}
\and {\it Volker Runde}\thanks{Research supported by NSERC under grant no.\ 227043-00.} \and {\it Nico Spronk}\thanks{The author wishes to thank his advisor, Brian Forrest, for making Oleg Aristov's visit to Waterloo possible.}}
\date{}
\begin{document}
\maketitle
\begin{abstract}
We investigate if, for a locally compact group $G$, the Fourier algebra $A(G)$ is biflat in the sense of quantized Banach homology. A central r\^ole in our investigation is played by the notion of an approximate indicator
of a closed subgroup of $G$: The Fourier algebra is operator biflat whenever the diagonal in $G \times G$ has an approximate indicator. Although we have been unable to settle the question of whether $A(G)$ is always operator biflat,
we show that, for $G = \SL(3,\comps)$, the diagonal in $G \times G$ fails to have an approximate indicator.
\end{abstract}
\begin{keywords}
locally compact groups, biflatness, Fourier algebra, (quantized) Banach homology, approximate indicator, Kazhdan's property $(T)$.
\end{keywords}
\begin{classification}
22D25 (primary), 22E10, 43A30, 46L07, 46L89, 46M18, 47L25, 47L50.
\end{classification}
\section*{Introduction}
In his seminal memoir \cite{Joh1}, B.\ E.\ Johnson proved that the amenable locally compact groups $G$ can be characterized by the vanishing of certain Hochschild cohomology groups of $L^1(G)$: This initiated the theory of
amenable Banach algebras. At about the same time, Banach homology, i.e.\ homological algebra with functional analytic overtones, was developed systematically by A.\ Ya.\
Helemski\u{\i}'s Moscow school (\cite{Hel}). One of the central notions in this theory is projectivity. A Banach algebra which is projective as a bimodule over itself is called biprojective. The biprojectivity
of $L^1(G)$ is equivalent to $G$ being compact (\cite[Theorem 51]{Hel2}). This shows that some important properties of locally compact groups $G$ are equivalent to certain homological properties of $L^1(G)$.
\par
The Fourier algebra $A(G)$ --- as introduced in \cite{Eym} --- can be viewed as the ``quantized'' counterpart of $L^1(G)$. In classical Banach homology, $A(G)$ fails to reflect the properties of $G$ in a satisfactory manner:
There are even compact groups $G$ for which $A(G)$ is not amenable (\cite{Joh2}). The reason for this failure lies in the fact that classical Banach homology ignores the canonical operator space structure, which $A(G)$
inherits as the predual of the group von Neumann algebra $\VN(G)$. The definition of an amenable Banach algebra, however, can easily be adapted to take operator space structures into account (\cite{Rua}): this leads to the notion of an
operator amenable, completely contractive Banach algebra. In \cite{Rua}, Z.-J.\ Ruan showed that in this quantized theory an analogue of Johnson's theorem holds for the Fourier algebra: $A(G)$ is operator amenable if and only if
$G$ is amenable. As in the classical theory, projectivity
plays a central r\^ole in quantized Banach homology. Dual to the classical situation, $A(G)$ is operator biprojective if and only if $G$ is discrete (\cite{Ari} and \cite{Woo}).
\par
Another important homological concept is that of biflatness. (A
Banach algebra is biflat if it is a flat Banach bimodule over
itself.) In general, biflatness is weaker than both amenability
and biprojectivity. Nevertheless, biflatness is of little
relevance in the study of group algebras: Since $L^1(G)$ always
has a bounded approximate identity, it is biflat if and only if
it is amenable (\cite[Theorem VII.2.20]{Hel}). This changes,
however, in the quantized setting: Since operator biprojectivity
implies operator biflatness, $A(G)$ is operator biflat whenever
$G$ is discrete. More generally, $A(G)$ is operator biflat for
every group $G$ such that $L^1(G)$ has a quasi-central bounded
approximate identity (\cite{RX}); this includes all amenable
groups (\cite[Theorem 3]{LR}), but also all $[\operatorname{SIN}]$-groups.
It is possible that every locally compact group has an operator biflat
Fourier algebra. One piece of (albeit circumstantial) evidence in
favor of this conjecture is the main result of \cite{Spr}: Every
Fourier algebra is operator weakly amenable, and, as in the
classical setting (\cite[Theorem 5.3.13]{Run}), operator
biflatness implies operator weak amenability.
\par
One goal of this paper is to systematically investigate whether
or not $A(G)$ is operator biflat for an arbitrary locally compact
group $G$. A central r\^ole in our investigation is played by the
notion of an approximate indicator: Roughly speaking, a closed
subgroup $H$ of a locally compact group $G$ has an approximate
indicator if its indicator function can be approximated in a
suitable way by functions from the Fourier--Stieltjes algebra. If the diagonal subgroup
of $G \times G$ has an approximate indicator, then $A(G)$ is
operator biflat. In particular, this is true if $G$ can be
continuously embedded into a group with a quasi-central bounded
approximate identity.
\par
The question of whether, for particular $G$ and $H$, an approximate indicator exists is of independent interest. We shall give both positive and negative results: Every neutral subgroup has an approximate indicator, but
$\SL(2,\comps)$ --- when canonically embedded into $\SL(3,\comps)$ --- fails to have one. Similarly, we shall see that the diagonal of $\SL(3,\comps) \times \SL(3,\comps)$ lacks an approximate indicator. This makes $\SL(3,\comps)$
a likely candidate for a group whose Fourier algebra is not operator biflat.
\subsection*{Acknowledgment}
While writing this paper, the authors have benefitted from conversations with and e-mails from the following people: Bachir Bekka, Brian Forrest, Eberhard Kaniuth, Tony Lau, Matthias Neufang, G\"unter Schlichting, George Willis, and Peter Wood.
Thanks to all of them!
\section{Flatness in the quantized setting}
The necessary background from (classical) Banach homology is covered in \cite{Hel} and, to a lesser extent, in \cite[Chapter 5]{Run}. Our reference for the theory of operator spaces is \cite{ER}, whose notation we adopt;
in particular, $\Tensor$ stands for the projective tensor product of operator spaces and not of Banach spaces.
\par
Since a lot of the development of Banach homology is categorical, many results carry over to the quantized,
i.e.\ operator space, context; for more details, see \cite{Rua}, \cite{RX}, \cite{WooPhD}, and \cite{Ari}, for example. We are therefore somewhat sketchy in our exposition here.
\begin{definition}
An algebra $\A$ which is also an operator space is called a {\it quantized Banach algebra\/} if the multiplication of $\A$ is a completely bounded bilinear map. If the multiplication is even completely contractive, $\A$
is called a {\it completely contractive Banach algebra\/}.
\end{definition}
\begin{remark}
For any quantized Banach algebra $\A$, multiplication induces a completely bounded linear map
\[
  \Delta \!: \A \Tensor \A \to \A, \quad a \tensor b \mapsto ab,
\]
the {\it diagonal map\/}.
\end{remark}
\begin{examples}
\item Let $\bar{\tensor}$ denote the $W^\ast$-tensor product. A {\it Hopf--von Neumann algebra\/} is a pair $(\M, \nabla)$, where $\M$ is a von Neumann algebra, and
$\nabla$ is a {\it co-multiplication\/}: a unital, $w^\ast$-continuous $^\ast$-monomorphism $\nabla \!: \M \to \M \bar{\tensor} \M$ which is co-associative, i.e.\ the diagram
\[
  \begin{CD}
  \M  @>\nabla>> \M \bar{\tensor} \M \\
  @V{\nabla}VV                                          @VV{\nabla \tensor \id_\M}V                     \\
  \M \bar{\tensor} \M @>>{\id_\M \tensor \nabla}> \M \bar{\tensor} \M \bar{\tensor} \M
  \end{CD}
\]
commutes. Let $\M_\ast$ denote the unique predual of $\M$. By \cite[Theorem 7.2.4]{ER}, we have $\M \bar{\tensor} \M \cong (\M_\ast \Tensor \M_\ast)^\ast$. Thus, $\nabla$ induces a
complete contraction $\nabla_\ast \!: \M_\ast \Tensor \M_\ast \to \M_\ast$ turning $\M_\ast$ into a completely contractive Banach algebra.
\item Let $G$ be a locally compact group, let $\wstar(G) := C^\ast(G)^{\ast\ast}$, and let $\omega \!: G \to \wstar(G)$ be the {\it universal representation\/} of $G$, i.e.\ for each (WOT-continuous and unitary) representation
$\pi$ of $G$ on a Hilbert space, there is unique $w^\ast$-continuous $^\ast$-homomorphism $\theta \!: \wstar(G) \to \pi(G)''$ such that $\pi = \theta \circ \omega$. Applying this universal property to the representation
\[
  G \to \wstar(G) \bar{\tensor} \wstar(G), \quad x \mapsto \omega(x) \tensor \omega(x)
\]
yields a co-multiplication $\nabla \!: \wstar(G) \to \wstar(G) \bar{\tensor} \wstar(G)$. Hence, $B(G) := \cstar(G)^\ast$, the {\it Fourier--Stieltjes\/} algebra of $G$, is a completely contractive Banach algebra.
\item The Fourier algebra $A(G)$ is a closed ideal of $B(G)$ (see \cite{Eym}), and thus also a completely contractive Banach algebra. (It is not hard to see that the operator space structure $A(G)$ inherited from
$B(G)$ coincides with the one it has the predual of $\VN(G)$.)
\end{examples}
\par
If $\A$ is a quantized Banach algebra, we call a left $\A$-bimodule $E$ a {\it quantized left $\A$-bimodule\/} if $E$ is an operator space such that the module operation is completely bounded; similarly, quantized right modules and bimodules
are defined. If $E$ is a quantized left $\A$-module, then $E^\ast$ with its right $\A$-module and operator space structures is a quantized right $\A$-module.
\par
In analogy with the classical situation, we call a short exact sequence
\begin{equation} \label{shortex}
  \{ 0 \} \to E_1 \to E_2 \to E_2/E_1 \to \{ 0 \}
\end{equation}
of quantized $\A$-modules (left, right, or bi-) {\it admissible\/} if $E_1$ is completely complemented in $E_2$ as an operator space, i.e.\ there is a completely bounded projection from $E_2$ onto $E_1$.
\par
Let $\A$ be a quantized Banach algebra, let $E$ be a quantized right $\A$-module, and let $F$ be a quantized left $\A$-module. Then $E \Tensor_\A F$ is defined as the quotient of $E \Tensor F$ modulo the closed linear span of the
set $\{ x \cdot a \tensor y - x \tensor a \cdot y : a \in A, \, x \in E, \, y \in F \}$.
\par
The notion of flatness for quantized modules is defined as in the classical setting:
\begin{definition} \label{flatdef}
Let $\A$ be a quantized Banach algebra. A quantized left $\A$-module $F$ is called {\it flat\/} if, for each short exact sequence (\ref{shortex}) of quantized right $\A$-modules, the complex
\[
  \{ 0 \} \to E_1 \Tensor_\A F \to E_2 \Tensor_\A F \to E_2/E_1 \Tensor_\A F \to \{ 0 \}
\]
is a short exact sequence of operator spaces.
\end{definition}
\par
Flatness for quantized right modules and bimodules is defined
analogously.
\par
Like flatness, the notion of injectivity (\cite[Definition 5.3.6]{Run}) translates to the quantized setting with the obvious modifications. As in the classical situation (\cite[Theorem 5.3.8]{Run}), we have (with a virtually identical proof):
\begin{theorem} \label{flatdual}
Let $\A$ be quantized Banach algebra, and let $E$ be a quantized left $\A$-module. Then the following are equivalent:
\begin{items}
\item $E$ is flat.
\item $E^\ast$ is an injective quantized right $\A$-module.
\end{items}
\end{theorem}
\par
With flatness for bimodules, the concept of a biflat Banach algebra carries over to the quantized situation:
\begin{definition}
A quantized Banach algebra $\A$ is called {\it operator biflat\/} if it is a flat quantized $\A$-bimodule.
\end{definition}
\par
The following characterization holds (with a proof analogous to
that of its classical counterpart from \cite{Hel}; for details, see \cite[Lemma 4.3.22 and Theorem 5.3.12]{Run}):
\begin{theorem} \label{biflatchar}
The following are equivalent for a quantized Banach algebra $\A$:
\begin{items}
\item $\A$ is operator biflat.
\item The adjoint $\Delta^\ast \!: \A^\ast \to (\A \Tensor \A)^\ast$ of the diagonal map has a completely bounded left inverse which is an $\A$-bimodule homomorphism.
\item There is a completely bounded $\A$-bimodule homomorphism $\rho \!: \A \to (\A \Tensor \A)^{\ast\ast}$ such that $\Delta^{\ast\ast} \circ \rho$ is the canonical embedding of $\A$ into $\A^{\ast\ast}$.
\end{items}
\end{theorem}
\par
The question which provided most of the motivation for this paper is whether or not, for a locally compact group $G$, its Fourier algebra $A(G)$ is operator biflat. For convenience, we define:
\begin{definition}
A locally compact group $G$ is called {\it biflat\/} if $A(G)$ is operator biflat.
\end{definition}
\par
We conclude this section with a hereditary property of biflat, locally compact groups:
\begin{proposition} \label{flatprop}
Let $G$ be a biflat, locally compact group, and let $H$ be a closed subgroup. Then $H$ is biflat, and $A(H)$ is a flat left quantized $\A$-module.
\end{proposition}
\begin{proof}
Let
\[
  I(H) := \{ f \in A(G) : f |_H \equiv 0 \}.
\]
Since $A(H) \cong A(G) / I(H)$ as operator spaces (\cite[Proposition 4.1]{Woo}), it is sufficient by (the quantized analogues of) \cite[Propositions 4 and 5]{Sel} to show that $A(G) \cdot I(H)$, i.e.\ the closed linear span of the set
$\{ fg : f \in A(G), \, g \in I(H) \}$ equals $I(H)$.
\par
Let $f \in I(H)$. Since $H$ is a set of synthesis for $A(G)$ by \cite[Theorem 2]{Her}, we can suppose that $f$ has compact support. Using the regularity of $A(G)$ (\cite[(3.2) Lemma]{Eym}), we find $g \in A(G)$ such that
$g |_{\supp \, f} \equiv 1$, so that $f = gf \in A(G) \cdot I(H)$. This proves the claim.
\end{proof}
\section{Approximate indicators and the operator biflatness of $A(G)$}
If $H$ is a closed subgroup of a locally compact group $G$, the indicator function of $H$ lies in $B(G)$ if and only if $H$ is open. In particular, each subgroup of a discrete group $G$ has
its indicator function in $B(G)$.
\par
We make the following definition:
\begin{definition} \label{appind}
Let $G$ be a locally compact group, and let $H$ be a closed
subgroup. A bounded net $( f_\alpha )_\alpha$ in $B(G)$ is called
an {\it approximate indicator\/} for $H$ if
\begin{alphitems}
\item $\lim_\alpha f (f_\alpha |_H) = f$ for all $f \in A(H)$;
\item $\lim_\alpha g f_\alpha = 0$ for all $g \in I(H)$.
\end{alphitems}
If $( f_\alpha )_\alpha$ is bounded by one, we speak of a {\it
contractive approximate indicator\/}, and if each function
$f_\alpha$ is positive definite, we call $( f_\alpha )_\alpha$ a
{\it positive definite approximate indicator\/}.
\end{definition}
\begin{remark}
In \cite{Joh3}, B.\ E.\ Johnson proved that a Banach algebra is amenable if and only if it has an approximate diagonal. An analogous statement holds in the quantized context (\cite{Rua}).
In terms of Definition \ref{appind}, \cite[Proposition 2.4 and Theorem 3.6]{Rua} assert that $G$ is amenable if and only if the diagonal subgroup
\[
  G_\Delta := \{ (x,x) : x \in G \}
\]
of $G \times G$ has an approximate indicator in $A(G \times G)$.
\end{remark}
\par
Recall that a closed subgroup $H$ of a locally compact group $G$
is called {\it neutral\/} if $e$ has a basis $\mathfrak U$ of
neighborhoods such that $U H = H U$ for all $U \in {\mathfrak U}$
(\cite[p.\ 96]{KL}); all normal subgroups are neutral, but the
same is true for every subgroup $H$ such that $e$ has a basis of neighborhoods invariant under
conjugation with elements of $H$.
\begin{example}
Let $G$ be a locally compact group, and let $H$ be a closed
neutral subgroup of $G$. For each compact subset $K$ of $G$ with
$K \cap H = \void$, \cite[Proposition 2.2]{KL} yields a
continuous, positive definite function $f_K$ on
$G$ with
\[
  f_K |_K \equiv 0 \qquad\text{and}\qquad f_K |_H \equiv 1.
\]
Let ${\mathbb K}$ be the collection of all compact subsets of $G$
which have empty intersection with $H$, ordered by set inclusion.
It is obvious that the net $( f_K )_{K \in {\mathbb K}}$
satisfies Definition \ref{appind}(a). Clearly, $\lim_K f g = 0$
holds for all $g \in A(G)$ with compact support disjoint from
$H$. Since $H$ is a set of synthesis for $A(G)$, Definition
\ref{appind}(b) holds as well. Hence, $( f_K )_{K \in {\mathbb
K}}$ is a (contractive, positive definite) approximate indicator
for $H$.
\end{example}
\begin{lemma} \label{biflatlem}
Let $G$ be a locally compact group, and let $H$ be a closed subgroup of $G$ which has an approximate indicator. Then there is
a completely bounded $A(G)$-module homomorphism $\rho \!: A(H) \to A(G)^{\ast\ast}$ such that $\Gamma_H^{\ast\ast} \circ \rho$ is
the canonical embedding of $A(H)$ into $A(H)^{\ast\ast}$, where $\Gamma_H \!: A(G)\to A(H)$ is the restriction map.
\end{lemma}
\begin{proof}
Let $( f_\alpha )_{\alpha \in {\mathbb A}}$ be an approximate indicator for $H$. For each $\alpha \in \mathbb A$,
define
\[
  \rho_\alpha \!: A(G)  \to B(G), \quad g \mapsto g f_\alpha.
\]
It is clear from this definition that each $\rho_\alpha$ is completely bounded with $\| \rho_\alpha \|_{\mathrm{cb}} \leq \| f_\alpha \|$ and an $A(G)$-module homomorphism. Since
$A(G) $ is an ideal in $B(G)$, each $\rho_\alpha$ attains its values in $A(G)$.
\par
Let $\cal U$ be an ultrafilter on $\mathbb A$ which dominates the
order filter, and define
\[
  \tilde{\rho} \!: A(G)  \to A(G )^{\ast\ast}, \quad  g \mapsto \text{$w^\ast$-}\lim_{\cal U} \rho_\alpha(g).
\]
Since $( f_\alpha )_{\alpha \in {\mathbb A}}$ is bounded, it is immediate that $\tilde{\rho}$ is well defined and
completely bounded; it is also clear that $\tilde{\rho}$ is an $A(G)$-module homomorphism. From Definition \ref{appind}(b) it
follows immediately that $\tilde{\rho}(g) = 0$ for all $g \in A(G)$ that vanish on $H$. Since $\Gamma_H \!: A(G) \to A(H)$ is a complete quotient map, it follows that
$\tilde{\rho}$ drops to a completely bounded $A(G)$-module homomorphism $\rho \!: A(H) \to A(G)^{\ast\ast}$.
\par
From Definition \ref{appind}(a), it is clear that $\Gamma_H^{\ast\ast} \circ \rho$ is the canonical embedding of
$A(H)$ into $A(H)^{\ast\ast}$.
\end{proof}
\par
The reason for our interest in approximate indicators stems from the following consequence of Lemma \ref{biflatlem}:
\begin{proposition} \label{biflatprop}
Let $G$ be a locally compact group such that $G_\Delta$ has an approximate indicator. Then $G$ is biflat.
\end{proposition}
\begin{proof}
By Lemma \ref{biflatlem} there is a completely bounded $A(G\times G)$-module homomorphism $\rho \!: A(G) \to A(G \times G)^{\ast\ast}$ such that $\Gamma_{G_\Delta}^{\ast\ast} \circ \rho$ is the
canonical embedding of $A(G)$ into $A(G)^{\ast\ast}$. Since $A(G \times G) \cong A(G) \Tensor A(G)$, and since  $\Gamma_{G_\Delta} = \Delta$, it follows from Theorem \ref{biflatchar} that $A(G)$ is operator biflat.
\end{proof}
\par
Recall that a locally compact group $G$ is called a
$[\SIN]$-group if $e$ has a basis of conjugation invariant
neighborhoods or, equivalently, if $L^1(G)$ has a bounded
approximate identity in its center. It is easy to see that $G$ is
a $[\SIN]$-group if and only if $G_\Delta$ is a neutral subgroup
of $G \times G$. Nevertheless, Proposition \ref{biflatprop}
allows us to establish the biflatness of locally compact groups
which are not $[\SIN]$-groups.
\par
We call a locally compact group $G$ a $[\operatorname{QSIN}]$-group --- $[\operatorname{QSIN}]$ standing for {\it quasi-$[\operatorname{SIN}]$\/} --- if $L^1(G)$ has a bounded approximate identity
$( e_\alpha )_{\alpha \in {\mathbb A}}$ such that
\[
  \delta_x \ast e_\alpha - e_\alpha \ast \delta_x \to 0 \qquad (x \in G).
\]
All $[\SIN]$-groups are trivially $[\operatorname{QSIN}]$-groups,
but so are all amenable groups (\cite[Theorem 3]{LR}); further results on
$[\operatorname{QSIN}]$-groups are contained in \cite{Sto}.
\par
Slightly extending \cite[Theorem 4.4]{RX}, we have:
\begin{theorem} \label{biflatthm}
Let $G$ be a locally compact group that can be continuously embedded into a $[\operatorname{QSIN}]$-group. Then $G$ is biflat.
\end{theorem}
\begin{proof}
We give an argument which avoids the Kac algebra machinery from \cite{RX} (see also \cite[Lemma~7.4.2]{Run}).
\par
Let $H$ be a $[\operatorname{QSIN}]$-group, and let $\theta \!: G \to H$ be an injective, continuous group homomorphism. By
\cite[Theorem 2]{LR}, we can find an  approximate identity $( e_\alpha )_{\alpha \in {\mathbb A}}$ for $L^1(H)$ with $e_\alpha \geq 0$ and $\| e_\alpha \|_1 = 1$ for all $\alpha \in \mathbb A$ such that
\[
  \|\delta_x \ast e_\alpha - e_\alpha \ast \delta_x \|_{L^1(H)} \to 0
\]
uniformly on compact subsets of $H$. Stokke's improved version of this result (\cite[Theorem 2.4]{Sto}) asserts
that $( e_\alpha )_{\alpha \in {\mathbb A}}$ can even be chosen with the supports tending to $\{ e_H \}$: for each neighborhood $U$
of $e_H$, there is $\beta \in \mathbb A$ such that $\supp \, e_\alpha\subset U$ for $\alpha \succcurlyeq \beta$.
\par
For each $\alpha \in \mathbb A$, let $\xi_\alpha := e_\alpha^\frac{1}{2}$. Let $\lambda$ and $\rho$ denote the
regular left and right representation, respectively, of $H$ on $L^2(H)$. Define
\[
  \mathbf{f}_\alpha(x,y) := \langle \lambda(\theta(x)) \rho(\theta(y)) \xi_\alpha, \xi_\alpha \rangle \qquad (x,y \in G, \, \alpha \in {\mathbb A}).
\]
\par
We claim that $( \mathbf{f}_\alpha )_{\alpha \in {\mathbb A}}$ is a contractive, positive definite approximate indicator for $G_\Delta$.
\par
Since
\begin{eqnarray*}
  |\mathbf{f}_\alpha(x,x)-1_{G\times G}|^2 & = & |\langle\lambda(\theta(x)) \rho(\theta(x)) \xi_\alpha, \xi_\alpha \rangle-\langle\xi_\alpha, \xi_\alpha \rangle|^2 \\
  & \leq & \|\lambda(\theta(x)) \rho(\theta(x)) \xi_\alpha- \xi_\alpha\|_{L^2(H)}^2 \\
  & = & \int_H \left( \xi_\alpha(\theta(x^{-1})\, y \, \theta(x) \Delta_H(\theta(x))^\frac{1}{2} - \xi_\alpha(y) \right)^2 dy, \\
  & & \qquad\text{where $\Delta_H$ is the modular function of $H$}, \\
  & \leq & \int_H | \xi_\alpha(\theta(x^{-1}) \, y \, \theta(x))^2 \Delta_H(\theta(x)) - \xi_\alpha(y)^2 | \, dy \\
  & = & \int_H | e_\alpha(\theta(x^{-1})\, y  \, \theta(x))\Delta_H(\theta(x)) - e_\alpha(y) | \, dy \\
  & = & \| \delta_{\theta(x)} \ast e_\alpha - e_\alpha \ast \delta_{\theta(x)} \|_{L^1(H)} \qquad (x \in G)
\end {eqnarray*}
(compare \cite[(3.8)]{Rua}), it follows that, $\mathbf{f}_\alpha(x,x) \to 1$ uniformly on compact subsets of $G$. By \cite[Theorem B$_2$]{GL}, this means that Definition \ref{appind}(a) holds.
\par
For Definition \ref{appind}(b), let  $\mathbf{g} \in I(G_\Delta)$. Since $G_\Delta$ is a set of synthesis for $A(G\times G)$ by \cite[Theorem 2]{Her}, we can suppose that
$\supp \, g \subset K$, where $K$ is a compact subset of $G\times G$ disjoint from $G_\Delta$. Let
\[
  \tilde{K} := \{ \theta(x^{-1}y): (x,y) \in K \}
\]
and
\[
  L := \{ x \in G: \text{there is $y \in G$ with $(x,y) \in K$} \}.
\]
Then $\tilde{K}$ and $L$ are compact sets, and since $\theta$ is injective, $e_H \notin \tilde{K}$ holds. Let $U$ be a
neighborhood of $e_H$ such that $U \cap Ux = \void$ for all $x \in \tilde{K}$. By \cite[4.5(iv)]{HR}, there is a neighborhood $V$ of $e_H$
such that $V \subset U$ and $\theta(x^{-1})V\theta(x)\subset U$ for all $x\in L$. Consequently,
\[
  \theta(x^{-1})V\theta(y)=\theta(x^{-1})V\theta(x)\theta(x^{-1}y)\subset U\theta(x^{-1}y)
  \qquad ((x,y)\in K).
\]
holds, so that $\theta(x^{-1})V\theta(y)\cap V=\void$ for all $(x,y) \in K$.
Let $\beta\in {\mathbb A}$ such that $\supp \, \xi_\alpha\subset V$ whenever $\alpha\succcurlyeq \beta$. If $(x,y)\in K$, then
$\supp \, \lambda(\theta(x)) \rho(\theta(y)) \xi_\alpha\subset  \theta(x)V\theta(y^{-1})$, and since $\theta(x)V\theta(y^{-1})\cap V=\void$, we have
\[
  \mathbf{f}_\alpha(x,y)= \langle\lambda(\theta(x)) \rho(\theta(y))
  \xi_\alpha, \xi_\alpha \rangle=0 \qquad ((x,y) \in K, \,
  \alpha \succcurlyeq \beta).
\]
Since $\mathbf{g}(x,y)=0$ for all $(x,y)\notin K$, this yields $\mathbf{g} \mathbf{f}_\alpha=0$ for $\alpha\succcurlyeq \beta$.
\end{proof}
\begin{remarks}
\item Theorem \ref{biflatthm} applies, in particular, to all locally
compact groups that can be continuously embedded into amenable
ones. Recall that a locally compact group is said to have {\it
Kazhdan's property $(T)$\/} if the trivial representation $1_G$
is an isolated point in $\hat{G}$. For example, $\SL(n,\free)$
has property $(T)$ for $n \geq 3$ and $\free = \reals$ or $\free
= \comps$ (\cite{dHV} and, in particular, \cite[Theorem 1.4.14]{BdHV}). Since $\SL(n,\free)$ has no
finite-dimensional unitary representations but the trivial one, it
follows from \cite[Corollary 7.1.10]{Zim} that $\SL(n,\free)$
does not continuously embed into an amenable, locally compact
group for $n \geq 3$. Hence, for such $G$, \cite[Theorem 3]{LR} and Theorem \ref{biflatthm} cannot not be used to establish the existence of an approximate indicator for $G_\Delta$. (In Theorem \ref{oleg4} below, we shall see
that for $G = \SL(3,\comps)$ even no such approximate indicator exists.)
\item We believe, but have been unable to prove, that there are indeed locally compact groups which continuously embed into $[\operatorname{QSIN}]$-groups without being $[\operatorname{QSIN}]$-groups themselves.
Potential candidates for such examples are \cite[Examples 1 and 2]{LR} (see \cite{LR1} for some of the technical details only sketched in \cite{LR}).
\end{remarks}
\section{The discretized Fourier--Stieltjes algebra}
In this section, we shall see that an approximate indicator --- even a contractive, positive definite one --- already exists if the indicator function of the subgroup under consideration can be approximated
by functions from the Fourier--Stieltjes algebra in a seemingly much weaker sense than required by Definition \ref{appind}(a) and (b).
\par
Let $G$ be a locally compact group, and let $G_d$ denote the same group, but equipped with the discrete topology. Recall that $B(G)$ is a closed subalgebra of $B(G_d)$ for any locally compact group which consists precisely
of the continuous functions in $B(G_d)$ (\cite[(2.24) Corollaire 1]{Eym}).
\begin{definition}
Let $G$ be a locally compact group. The {\it discretized Fourier--Stieltjes algebra\/} $B_d(G)$ of $G$ is the $w^\ast$-closure of $B(G)$ in $B(G_d)$.
\end{definition}
\begin{remarks}
\item It is immediate that $B_d(G)$ is a $w^\ast$-closed subalgebra of $B(G_d)$ whose predual is the $\cstar$-subalgebra of $\wstar(G)$ generated by the set $\{ \omega(x) : x \in G \}$, which we denote by $C^\ast_d(G)$.
(It is easy to see
that $B_d(G)$ is a dual Banach algebra in the sense of \cite{Run2}.)
\item Since the embedding of $C^\ast_d(G)$ into $\wstar(G)$ is an injective $^\ast$-homomorphism and thus an isometry, the inclusion of $B(G)$ into $B_d(G)$ extends to a $w^\ast$-continuous quotient map from $B(G)^{\ast\ast}$
onto $B_d(G)$.
\item It may well be that $B_d(G) \subsetneq B(G_d)$: For a connected Lie group $G$, the equality $B_d(G) = B(G_d)$ holds if and only if $G$ is solvable (\cite{BV}).
\item Trivially, $B_d(G) = B(G_d)$ holds if $G$ is discrete. The same is true if $G_d$ is amenable (\cite[Corollary 1.5]{BLS}).
\item In \cite{BKLS}, it is conjectured that $B_d(G) = B(G_d)$ is true if and only if $G$ contains an open subgroup which is amenable as a discrete group.
\end{remarks}
\par
Let $G$ be a locally compact group and let $H$ be a closed subgroup of $G$ which has an approximate indicator. It is straightforward to see that then $\chi_H$ belongs to $B_d(G)$. In the remainder of this section, we shall see that the
converse holds as well.
\par
We first introduce some notation.
\par
Let $G$ be a locally compact group, and let $S \subset G$ be any subset. We define
\[
  I(B(G),S) := \{ f \in B(G) : f|_S \equiv 0 \} \qquad\text{and}\qquad I(B_d(G),S) := \{ f \in B_d(G) : f|_S \equiv 0 \}.
\]
Moreover, let $C^\ast_d[S]$ and $\wstar[S]$ denote the norm closed and the $w^\ast$-closed linear span of $\{ \omega(x) : x \in S \}$ in $\wstar(G)$, respectively. It follows from the bipolar theorem that
\begin{equation} \label{polars1}
  I(B_d(G),S))^\circ = C^\ast_d[S] \qquad\text{and}\qquad I(B(G),S))^\circ = \wstar[S] \cap C^\ast_d(G),
\end{equation}
where the polar is taken in $C^\ast_d(G)$.
\par
We proceed by proving a series of lemmas (plus one corollary and one proposition). All polars are taken with respect to the canonical duality between $C^\ast_d(G)$ and $B_d(G)$.
\begin{lemma} \label{discrlem1}
Let $G$ be a locally compact group, let $U \subset G$ be open, and let $K \subset U$ be compact. Then $I(B_d(G),U) \subset I(B(G),K)^{\circ\circ}$ holds.
\end{lemma}
\begin{proof}
Let $f \in I(B_d(G),U)$. Using the regularity of the Fourier
algebra (\cite[(3.2) Lemme]{Eym}), we find $g \in A(G)$ such that
$g |_K \equiv 1$ and $g |_{G \setminus U} \equiv 0$. Let $(
f_\alpha )_\alpha$ be a net in $B(G)$ such that $f_\alpha \to f$
in the $w^\ast$-topology on $B_d(G)$. It follows that $f_\alpha
(g-1) \to f(g-1) = f$ in the $w^\ast$-topology. Hence, $f$ lies
in the $w^\ast$-closure of $I(B(G),K)$ in $B_d(G)$. The bipolar
theorem then yields the claim.
\end{proof}
\par
From the polar point of view, Lemma \ref{discrlem1} reads as:
\begin{corollary} \label{discrcor1}
Let $G$ be a locally compact group, let $U \subset G$ be open, and let $K \subset U$ be compact. Then $\wstar[K] \cap C^\ast_d(G) \subset C^\ast_d[U]$ holds.
\end{corollary}
\begin{lemma} \label{discrlem2}
Let $G$ be a locally compact group, let $H$ be a closed subgroup such that $\chi_H \in B_d(G)$, and let ${\cal L}_H$ denote the collection of all compact subsets of $G$ which have empty intersection with $H$. Then
\begin{equation} \label{discreq1}
  \bigcap_{L \in {\cal L}_H} C^\ast_d[G \setminus L] = C^\ast_d[H]
\end{equation}
holds.
\end{lemma}
\begin{proof}
First note that trivially
\[
  \bigcap_{L \in {\cal L}_H} I(B_d(G),L) = I(B_d(G),G\setminus H)
\]
holds. Let
\[
  J := \bigcup_{L \in {\cal L}_H} I(B_d(G),G \setminus L).
\]
It is immediate that $J$ is an ideal of $B_d(G)$. \par Let $L'
\in {\cal L}_H$. Choose $L \in {\cal L}_H$ with $L'$ in the
interior of $L$. The regularity of $A(G)$ yields $f \in A(G)$
with $f |_{L'} \equiv 1$ and $f |_{G \setminus L} \equiv 0$.
Hence, $f \in J$ and $1-f \in I(B_d(G),L')$. It follows that $1
\in J + I(B_d(G),L')$ and thus $J + I(B_d(G),L') = B_d(G)$. In
terms of polars, this means that
\[
  \bigcap_{L \in {\cal L}_H} C^\ast_d[G \setminus L] \cap C^\ast_d[L'] = \{ 0 \}.
\]
Since $L' \in {\cal L}_H$ was arbitrary, we obtain
\[
  \bigcap_{L \in {\cal L}_H} C^\ast_d[G \setminus L] \cap \bigcup_{L' \in {\cal L}_H} C^\ast_d[L'] = \{ 0 \}
\]
and thus --- taking polars again --- that $J + I(B_d(G),G \setminus H)$ is $w^\ast$-dense in $B_d(G)$. Since $\chi_H \in B_d(G)$, we have $\chi_{G\setminus H} \in B_d(G)$ as well. Multiplication in $B_d(G)$ is separately continuous;
we therefore obtain that
\[
  \chi_{G \setminus H}(J + I(B_d(G),G \setminus H)) = \chi_{G \setminus H} J = J
\]
is $w^\ast$-dense in $\chi_{G\setminus H} B_d(G) = I(B_d(G),H)$. Taking polars for one last time (in this proof at least) yields (\ref{discreq1}).
\end{proof}
\begin{lemma} \label{discrlem3}
Let $G$ be a locally compact group, let $H$ be a closed subgroup of $G$ such that $\chi_H \in B_d(G)$, and let $K \subset H$ be compact. Then $\wstar[K] \cap C^\ast_d(G) \subset C^\ast_d[H]$ holds.
\end{lemma}
\begin{proof}
Let ${\cal L}_H$ be as in Lemma \ref{discrlem2}. For any $L \in {\cal L}_H$, Corollary \ref{discrcor1} yields
\[
  \wstar[K] \cap C^\ast_d(G) \subset C^\ast_d[G \setminus L].
\]
Since $L \in {\cal L}_H$ was arbitrary, this in turn implies
\[
  \wstar[K] \cap C^\ast_d(G) \subset \bigcap_{L \in {\cal L}_H} C^\ast_d[G \setminus L] = C^\ast_d[H]
\]
by Lemma \ref{discrlem2}.
\end{proof}
\par
For any normed space $E$, we denote its closed unit ball by $B_1[E]$.
\begin{proposition} \label{discrprop}
Let $G$ be a locally compact group, let $H$ be a closed subgroup of $G$ such that $\chi_H \in B_d(G)$, let $K \subset H$ and $L \subset G \setminus H$ be compact
and let $\epsilon >0$. Then there is $f \in B(G)$ with $\| f \| \leq 1$ such that
\[
  f |_L \equiv 0 \qquad\text{and}\qquad |f(x) -1 | < \epsilon \quad( x \in K).
\]
\end{proposition}
\begin{proof}
Let
\[
  S := \{ f \in B(G) : \text{$\| f \| \leq 1$, $f |_L \equiv 0$, and $f$ is constant on a neighborhood of $K$} \}.
\]
Then
\[
  S = S_1 \cap S_2 \cap B_1[B(G)]
\]
holds, where $S_1 = I(B(G),L)$ and
\[
  S_2 = I(B(G),K) + \comps.
\]
It follows that
\begin{eqnarray*}
  S^\circ & = & \varcl{S_1^\circ + S_2^\circ +  B_1[B(G)]^\circ}^{\| \cdot \|} \\
          & = & \varcl{S_1^\circ + S_2^\circ +  B_1[C^\ast_d(G)]}^{\| \cdot \|}
\end{eqnarray*}
and therefore
\begin{eqnarray}
  S^{\circ\circ} & = & S_1^{\circ\circ} \cap S_2^{\circ\circ} \cap B_1[C^\ast_d(G)]^\circ \nonumber \\
  & = & S_1^{\circ\circ} \cap S_2^{\circ\circ} \cap B_1[B_d(G)]. \label{polars4}
\end{eqnarray}
By Lemma \ref{discrlem1}, we have $S_1^{\circ\circ} \supset
I(B_d(G),G\setminus L)$, and Lemma \ref{discrlem3} yields
$S_2^{\circ\circ} \supset I(B_d(G),H) + \comps$. Hence, $\chi_H$
lies in the right hand side of (\ref{polars4}) and thus in
$S^{\circ\circ}$.
\par
By the bipolar theorem, $\chi_H$ therefore lies in the $w^\ast$-closure of $S$ in $B_d(G)$. Fix $x_0 \in K$. Then there is $f \in S$ --- which means, in particular, that $\| f \| \leq 1$ and $f |_L \equiv 0$ --- with $| f(x_0) - 1   | < \epsilon$.
Since $f$ is constant on (some neighborhood of) $K$, we have in fact that $| f(x) - 1 | < \epsilon$ for all $x \in K$.
\end{proof}
\begin{theorem} \label{discrthm}
Let $G$ be a locally compact group and let $H$ be a closed subgroup of $G$. Then the following are equivalent:
\begin{items}
\item $\chi_H \in B_d(G)$.
\item There is a contractive approximate indicator for $H$.
\item There is a contractive, positive definite approximate indicator for $H$.
\item There is an approximate indicator for $H$.
\end{items}
\end{theorem}
\begin{proof}
(i) $\Longrightarrow$ (ii): By Proposition \ref{discrprop}, there is a net $( f_\alpha )_\alpha$ in $B(G)$ bounded by one that satisfies Definition \ref{appind}(b). Since $f_\alpha |_H \to 1$ uniformly on compact subsets of
$H$, it follows from \cite[Theorem B$_2$]{GL} that Definition \ref{appind}(a) is also satisfied.
\par
(ii) $\Longrightarrow$ (iii): Let $( f_\alpha )_{\alpha \in {\mathbb A}}$ be a contractive approximate indicator for $H$. It is clear that then $( f^\ast_\alpha )_\alpha$ is also a contractive approximate indicator
for $H$ as is $\left( \frac{1}{2}(f_\alpha + f^\ast_\alpha) \right)_{\alpha \in {\mathbb A}}$. We may thus suppose that $( f_\alpha )_\alpha$ consists of self-adjoint elements of $B(G)$. For each $\alpha \in \mathbb A$, let
$f_\alpha^+$ and $f_\alpha^-$ be the positive and negative part of $f_\alpha$, respectively, i.e.\ $f_\alpha^+$ and $f_\alpha^-$ are positive definite such that
\begin{equation} \label{lebesgue}
  f_\alpha = f_\alpha^+ - f_\alpha^- \qquad\text{and}\qquad \| f_\alpha \| = \| f_\alpha^+ \| + \| f_\alpha^- \|.
\end{equation}
Passing to subnets we can suppose that $( f^\pm_\alpha )_\alpha$ have $w^\ast$-limits $f^\pm$ in $B_d(G)$. From (\ref{lebesgue}) it then follows, in particular, that
\[
  1 = \lim_\alpha f_\alpha(e) = \lim_\alpha f^+_\alpha(e) - \lim_\alpha f^-(e) = f^+(e) - f^-(e).
\]
Since $\| f^+_\alpha \| \leq 1$ for all $\alpha \in \mathbb A$, we have $f^+(e) = \| f^+ \| \leq 1$. Since $f^-(e) \geq 0$, this necessitates that $f^-(e) = 0$ and thus $f^+(e) = 1$. Therefore,
\[
  \lim_\alpha \| f_\alpha^+ \| = \lim_\alpha f^+_\alpha(e) = f^+(e) = 1
\]
holds. Since $\| f_\alpha \| \leq 1$ for all $\alpha \in \mathbb A$, this means that $\lim_\alpha \| f^-_\alpha \| = 0$. Consequently, $( f^+_\alpha )_\alpha$ is an approximate indicator for $H$.
\par
(iii) $\Longrightarrow$ (iv) is trivial, and (iv) $\Longrightarrow$ (i) is straightforward, as was previously observed.
\end{proof}
\begin{remark}
It is easy to see from the proofs of Lemma \ref{biflatlem} and Proposition \ref{biflatprop} that the $\mathrm{cb}$-norm of the corresponding
splitting morphism $\rho \!: A(G) \to A(G\times G)^{\ast\ast}$ is less than or equal to any bound for the approximate indicator. It therefore follows from Theorem \ref{discrthm}
that the existence of approximate indicator already implies that we can find a splitting morphism $\rho \!: A(G) \to A(G\times G)^{\ast\ast}$ with $\|\rho\|_{\mathrm{cb}} \leq 1$.
\end{remark}
\par
The following example shows that, in general, Proposition \ref{discrprop} cannot be improved to guarantee the existence of $f \in B(G)$ with $\| f \| \leq 1$, $f |_L \equiv 0$ such that $|f(x) - 1| < \epsilon$ for all $x \in H$:
\begin{example}
Let $G$ be the $ax+b$ group, i.e.\
\[
  G = \{ (a,b) : a, b \in \reals, a > 0 \}
\]
with multiplication
\[
  (a_1,b_1)(a_2,b_2) := (a_1 a_2, b_1 + a_1 b_2) \qquad ((a_1, b_1), (a_2, b_2) \in G),
\]
and let $H = \{ (a,0) : a > 0 \}$. It will follow from Proposition \ref{oleg1} below that $H$ has an approximate indicator, say $( f_\alpha )_{\alpha \in {\mathbb A}}$, so that $\chi_H \in B_d(G)$. Assume now that
$f_\alpha |_H \to 1$ uniformly on $H$; suppose without loss of generality that $( f_\alpha )_\alpha$ is positive definite. It can be shown that $f_\alpha$ converges to $0$ uniformly on
$H (a,b) H$ for all $(a,b) \in G \setminus H$. Hence, for sufficiently large $\alpha \in \mathbb A$, the function $f_\alpha$ is arbitrarily small on $H(a,b)H$ and
arbitrarily close to $1$ on $H$. As pointed out in \cite[Example 1.3(i)]{KL}, the closure of $H(1,b)H$ for any $b > 0$ has non-empty intersection with $H$. This yields a contradiction due to the continuity of $f_\alpha$
for $\alpha \in \mathbb A$.
\end{example}
\section{Existence of approximate indicators}
Apart from the problem of whether or not every locally compact group is biflat, the question of whether a particular closed subgroup $H$ of a locally compact group $G$ has an approximate indicator is an intriguing question by itself.
\par
In this section, we shall give necessary conditions for the existence of an approximate indicator, but also encounter examples, where no approximate indicator can exist.
\par
Let $G$ be a locally compact groups, and let $\lambda$ denote the left regular representation of $G$ on $L^2(G)$. For any closed subgroup $H$ of $G$, let $\VN[H]$ be the $w^\ast$-closed, linear span of $\lambda(H)$ in $\VN(G)$.
It is well known that $\VN[H] \cong \VN(H)$.
\begin{proposition} \label{oleg1}
Let $G$ be a locally compact group, and let $H$ be a closed subgroup of $G$ such that there is a projection ${\cal E} \!: \VN(G) \to \VN[H]$ which is an $A(G)$-bimodule homomorphism. Then the following are equivalent:
\begin{items}
\item There is an approximate indicator for $H$ in $A(G)$.
\item $H$ is amenable.
\end{items}
\end{proposition}
\begin{proof}
(i) $\Longrightarrow$ (ii): Suppose that $H$ has an approximate indicator $( f_\alpha )_\alpha$ in $A(G)$. Then Definition \ref{appind}(a) yields that $( f_\alpha |_H )_\alpha$ is a bounded approximate identity for $A(H)$, so that
$H$ is amenable by Leptin's theorem (\cite[Theorem 7.1.3]{Run}).
\par
(ii) $\Longrightarrow$ (i). Suppose that $H$ is amenable. Leptin's theorem yields an approximate identity $( e_\alpha )_\alpha$ for $A(H)$. Let $E \in A(H)^{\ast\ast}$ be a $w^\ast$-accumulation point of $( e_\alpha )_\alpha$, and
let $F := {\cal E}^\ast E \in A(G)^{\ast\ast}$. It follows that
\[
  f \cdot F = 0 \quad (f \in I(H)) \qquad\text{and}\qquad F |_{\VN[H]} = E.
\]
Choosing a net in $A(G)$ that converges to $F$ in the $w^\ast$-topology and passing to convex combinations, we obtain
an approximate indicator for $H$ (compare \cite{Joh3}, where the equivalence of the existence of an approximate diagonal and the existence of a virtual diagonal is proved).
\end{proof}
\begin{remarks}
\item A module homomorphism $\cal E$ as in Proposition \ref{oleg1} exists in any of the following situations:
\begin{enumerate}
\item Suppose that $G$ has the $H$-separation property introduced and studied in \cite{KL}, i.e.\ for each $x \in G \setminus H$, there is a positive definite, continuous function $f$ on $G$ such that
\[
  f |_H \equiv 1 \qquad\text{and}\qquad f(x) \neq 1.
\]
By \cite[Proposition 3.1]{KL}, there is a norm one projection ${\cal E} \!: \VN(G) \to \VN[H]$ which is an $A(G)$-module homomorphism. By \cite[Theorem 5.1.5]{Li}, $\cal E$ is completely positive and thus completely bounded.
\item Suppose that $G$ is biflat, and that $\VN(H) \cong \VN[H]$ is injective (as a von Neumann algebra), e.g.\ if $H$ is amenable or connected (\cite[(1.31)]{Pat}). From the definition of injectivity for von Neumann algebras,
it is immediate that there is a (necessarily completely bounded)
norm one projection ${\cal E}_0 \!: \VN(G) \to \VN[H]$, i.e.\ the
short exact sequence
\[
  \{0 \} \to \VN[H] \to \VN(G) \to \VN(G) / \VN[H] \to \{ 0 \}
\]
of quantized left $A(G)$-modules is admissible. Since $A(H)$ is a flat quantized $A(G)$-module by Proposition \ref{flatprop}, $\VN(H)$ is injective (as a quantized $A(G)$-module) by Theorem \ref{flatdual}. It follows that the identity on
$\VN[H]$ extends to a completely bounded module homomorphism ${\cal E} \!: \VN(G) \to \VN(H)$.
\item Suppose that $H$ has an approximate indicator (not necessarily in $A(G)$), and let $\rho \!: A(H) \to A(G)^{\ast\ast}$ be the $A(G)$-module homomorphism that exists according to Lemma \ref{biflatlem}. Letting
${\cal E} := \rho^\ast |_{\VN(G)}$, we obtain a completely bounded $A(G)$-module homomorphism as in Proposition \ref{oleg1}. As a consequence of Theorem \ref{discrthm} and the remark following it, we can even make sure that
$\| {\cal E} \|_{\mathrm{cb}} = 1$.
\end{enumerate}
\item For general $G$ and $H$, there need not be a norm one projection ${\cal E} \!: \VN(G) \to \VN[H]$, let alone one that is an $A(G)$-module homomorphism. This can be seen as follows. The group $G = \SL(2,\reals)$ is connected,
so that $\VN(G)$ is an injective von Neumann algebra. Hence, there is a norm one projection from ${\cal B}(L^2(G))$ onto $\VN(G)$. It is well known that $G$ contains $\free_2$, the free group in two generators as a closed subgroup
(\cite[(3.2) Propsosition]{Pat}). Hence, if there were a norm one projection from $\VN(G)$ onto $\VN[\free_2]$, there would be such a projection from ${\cal B}(L^2(G))$ onto $\VN[\free_2]$, so that $\VN(\free_2)$ would be an injective
von Neumann algebra; this, in turn, would imply the amenability of $\free_2$ by \cite[Theorem 4.4.13]{Run}, which is wrong. In particular, $\free_2$ does not have an approximate indicator in $B(\operatorname{SL}(2,\reals))$.
\end{remarks}
\par
The question of whether, for general $G$ and $H$, a (not necessarily completely) bounded $A(G)$-module projection from $\VN(G)$ onto $\VN[H]$ exists seems to be open. Therefore, the following theorem on the existence of approximate indicators is
of interest:
\begin{theorem} \label{oleg2}
Let $G$ be a unimodular, locally compact group, let $H$ and $L$
be closed subgroups such that $H$ is amenable, $L$ is unimodular,
and
\[
  m \!: L \times H \to G, \quad (x,y) \mapsto xy
\]
is a homeomorphism. Then there is a contractive, positive definite approximate indicator for $H$ in $A(G)$.
\end{theorem}
\begin{proof}
Let $\mathbb A$ be the collection of all triples $(K,C,\epsilon)$, where $K \subset H$ and $C \subset G \setminus H$ are compact and $\epsilon > 0$. There is a natural order on $\mathbb A$, namely
\[
  (K_1,C_1,\epsilon_1) \preccurlyeq (K_2,C_2,\epsilon_2) \defiff \text{$K_1 \subset K_2$, $C_1 \subset C_2$, and $\epsilon_1 \geq \epsilon_2$}.
\]
\par
Let $\alpha = (K,C,\epsilon) \in \mathbb A$, and let $C_H$ and
$C_L$ be the projections of $m^{-1}(C)$ onto $H$ and $L$,
respectively. Since $H$ is amenable, it satisfies the F{\o}lner
condition (\cite[(4.10) Theorem]{Pat}), i.e.\ there is a compact
set $U_\alpha$ in $H$ such that $\| L_x \phi_\alpha - \phi_\alpha
\|_{L^1(H)} < \epsilon$ for all $x \in K$, where $\phi_\alpha$ is
the $L^1(H)$-normalized indicator function of $U_\alpha$ and
$L_x$ denotes left translation by $x$.
\par
Since $C_L C_H U_\alpha \cap U_\alpha = \void$, there is a compact
neighborhood $V_\alpha$ of $e \in L$ such that
\begin{equation} \label{void}
  C_L C_H U_\alpha V_\alpha \cap U_\alpha V_\alpha = \void
\end{equation}
Let $\xi_\alpha$ be the $L^2(G)$-normalized indicator function of $U_\alpha V_\alpha$, and let
\[
  f_\alpha(x) := \langle \lambda(x) \xi_\alpha, \xi_\alpha \rangle \qquad (x \in G).
\]
\par
We claim that $( f_\alpha )_\alpha$ is an approximate indicator for $H$.
\par
It follows immediately from (\ref{void}) that Definition \ref{appind}(b) is satisfied.
\par
Since $G$ and $L$ are both unimodular, Haar measure on $G$ can be identified with the product of Haar measures on $H$ and $L$, respectively (\cite[Chap\^{\i}tre VII, {\S}2, Proposition 13]{Bou}). For $\alpha = (K,C, \epsilon )$, we thus have
\begin{eqnarray*}
 | \langle L_x \xi_\alpha,\xi_\alpha \rangle - \langle\xi_\alpha, \xi_\alpha \rangle|^2 & \leq & \| L_x \xi_\alpha - \xi_\alpha \|_{L^2(G)}^2 \\
 & \leq & \int_G | \xi_\alpha(x^{-1}y)^{2}- \xi_\alpha(y)^2 | \, dy \\
 & = & \| L_x \xi_\alpha^2 - \xi_\alpha^2 \|_{L^1(G)} \\
 & = & \| L_x \phi_\alpha - \phi_\alpha \|_{L^1(H)} \\
 & < & \epsilon \qquad (x \in K),
\end {eqnarray*}
which entails that $f_\alpha |_H \to 1$ uniformly on compact subsets of $H$. Again, \cite[Theorem B$_2$]{GL} establishes Definition \ref{appind}(a).
\end{proof}
\begin{example}
Let $G = \operatorname{SL}(n,\reals)$ ($n\ge 2$), and let $H$ be
the subgroup of $G$ consisting of all upper triangular matrices
with positive diagonal elements. Then Theorem \ref{oleg2} is
applicable (with $L = \operatorname{SO}(n)$), so that there is an
approximate indicator for $H$.
\end{example}
\par
We have already remarked that $\free_2$, if embedded into $\operatorname{SL}(2,\reals)$ as a closed subgroup, cannot have an approximate indicator. We now turn to proving further (and more interesting) non-existence results.
\par
In what follows we shall always consider $\SL(2,\comps)$ as subgroup of $\SL(2,\comps)$ via
\[
  \SL(2,\comps) \to \SL(3,\comps), \quad \left[ \begin{array}{cc} a & b \\ c & d \end{array} \right] \mapsto \left[ \begin{array}{ccc} a & b & 0 \\ c & d & 0 \\ 0 & 0 & 1 \end{array} \right]
\]
\begin{proposition} \label{oleg3}
Let $G = \SL(3,\comps)$, let $H = \SL(2,\comps)$, and suppose that there is a bounded net $( f_\alpha )_{\alpha \in {\mathbb A}}$ of continuous, positive definite functions such that $f_\alpha |_H \to 1_H$ uniformly on compact subsets if $H$.
Then $f_\alpha \to 1_G$ in the norm topology.
\end{proposition}
\begin{proof}
For each $\alpha \in \mathbb A$, let $\pi_\alpha$ be a unitary
representation of $G$ of which $f_\alpha$ is a coefficient
function. Since $G$ has Kazhdan's property $(T)$, we have, for
each $\alpha \in \mathbb A$, a decomposition $\pi_\alpha =
\pi_{\alpha,1} \oplus \pi_{\alpha,2}$, where $\pi_{\alpha,1}
\prec  \hat{G} \setminus \{ 1_G \}$ (with $\prec$ standing for
weak containment) and $\pi_{\alpha,2} \prec 1_G$. Consequently,
we have a decomposition
\[
  f_\alpha(x) = f_{\alpha,1}(x) + t_\alpha \qquad (x \in G, \, \alpha \in {\mathbb A}),
\]
where $f_{\alpha,1}$ is a positive definite function associated with $\pi_{\alpha,1}$ and $t_\alpha \geq 0$.
\par
Fell's theorem (\cite[Theorem 6.1]{Fel}; see also \cite[Remark 1.13]{BLS}) asserts that $\pi |_H \prec \lambda_H$, where $\lambda_H$ is the regular left representation of $H$, for all $\pi \in \hat{G}  \setminus \{ 1_G \}$. It follows that
\[
  \pi_{\alpha,1} |_H \prec \{ \pi |_H : \pi \in \hat{G} \setminus \{ 1_G \} \} \prec \lambda_H \qquad (\alpha \in {\mathbb A}).
\]
\par
Since the net $( t_\alpha )_\alpha$ is bounded, we may replace it
by a subnet and suppose that it converges to some $t \geq 0$. If
$t \neq 1$, then $( f_{\alpha,1} |_H)_\alpha$ converges to a
non-zero constant function uniformly on compact subsets, so that
$1_H \prec \lambda_H$. But this means that $H$ is amenable
(\cite[p.\ 144]{Pat}), which is clearly false. Hence, $t =1$ must
hold. This, in turn, implies that $\| f_{\alpha,1} \| =
f_{\alpha,1}(e) \to 0$. It follows that $f_\alpha \to 1_G$ in the
norm topology of $B(G)$.
\end{proof}
\begin{corollary}
Let $G = \SL(3,\comps)$, and let $H = \SL(2,\comps)$. Then $H$ does not have an approximate indicator.
\end{corollary}
\begin{proof}
Assume towards a contradiction that the claim is false. By Theorem \ref{discrthm}, we can then suppose that $H$ has a positive definite approximate indicator, say $( f_\alpha )_\alpha$. From Definition \ref{appind}(a),
it follows that $f_\alpha |_H \to 1_H$ uniformly on compact subsets of $H$. By Proposition \ref{oleg3}, however, this violates Definition \ref{appind}(b).
\end{proof}
\par
We conclude this paper with another negative result, which casts doubt on whether $\SL(3,\comps)$ is a biflat group:
\begin{theorem} \label{oleg4}
Let $G = \SL(3,\comps)$. Then $G_\Delta$ does not have an approximate indicator.
\end{theorem}
\begin{proof}
First note that $G$ is a type I group (\cite{Dix}), so that $\cstar(G)$ is of type I and, in particular, nuclear. We therefore have that
\[
  \cstar(G \times G) = \cstar(G) \ttensor_{\max} \cstar(G) = \cstar(G) \ttensor_{\min} \cstar(G),
\]
where $\ttensor_{\max}$ and $\ttensor_{\min}$ denote the maximal and the minimal $\cstar$-tensor product, respectively. It follows (\cite[Lemma 13.5.6]{Li}) that
\[
  \widehat{G \times G} \cong \hat{G} \times \hat{G}.
\]
\par
Assume that $G_\Delta$ has an approximate indicator $( {\mathbf f}_\alpha )_{\alpha \in {\mathbb A}}$, which we can suppose to be positive definite by Theorem \ref{discrthm}. For each $\alpha \in \mathbb A$, let $\pi_\alpha$ be a unitary
representation of $G \times G$ with which ${\mathbf f}_\alpha$ is associated. Let
\[
  S_1 := \{ 1_G \times 1_G \}, \quad S_2 := (\hat{G} \setminus \{ 1_G\} ) \times \{ 1_G \}, \quad S_3 := \{ 1_G \} \times (\hat{G} \setminus \{ 1_G\})
\]
and
\[
  S_4 :=  (\hat{G} \setminus \{ 1_G \} ) \times (\hat{G} \setminus \{ 1_G \}).
\]
For each $\alpha \in {\mathbb A}$, there are unitary representations $\pi_{\alpha,1}, \ldots, \pi_{\alpha,4}$ of $G \times G$ such that
\[
  \pi_\alpha = \pi_{\alpha,1} \oplus \cdots \oplus \pi_{\alpha,4} \qquad\text{and}\qquad \pi_{\alpha,j} \prec S_j \quad (j=1,2,3,4).
\]
For notational simplicity, we write $G$ instead of $G_\Delta$ for the remainder of this proof. It is clear that
\[
  \pi_{\alpha,1} |_G \prec 1_G \qquad\text{and}\qquad \pi_{\alpha,j} |_G \prec \hat{G} \setminus \{ 1_G \} \quad (j=2,3)
\]
for all $\alpha \in \mathbb A$.
\par
We claim that $\pi_{\alpha,4} |_G \prec \hat{G} \setminus \{ 1_G \}$ holds as well for all $\alpha \in \mathbb A$. Assume that $\pi_{\alpha,4} |_G \not\prec \hat{G} \setminus \{ 1_G \}$ for some $\alpha \in \mathbb A$. Since
$G$ has property $(T)$, it is easy to see that then $1_G \prec \pi_{\alpha,4} |_G$ holds. Let $\Pi$ be the sum of all representations of $G$ in $\hat{G} \setminus \{ 1_G \}$. It follows that $\overline{\Pi} = \Pi$ and that $\pi_{\alpha,4}
\prec \Pi \times \Pi = \Pi \times \overline{\Pi}$. Consequently,
\[
  1_G \prec \pi_{\alpha,4} |_G \prec \Pi \tensor \overline{\Pi}
\]
holds. From \cite[Theorem 5.1]{Bek}, it follows that $\Pi$ is an amenable representation of $G$ in the sense of \cite[Definition 1.1]{Bek}. Since $G$ has Kazhdan's property $(T)$, \cite[Corollary 5.9]{Bek} implies that $\Pi$ has a finite-dimensional,
unitary subrepresentation. Since $1_G$ is the only such representation of $G$, this means that $1_G \prec \Pi \prec \hat{G} \setminus \{ 1_G \}$, so that we have reached a contradiction.
\par
As in the proof of Proposition \ref{oleg3}, we thus obtain a decomposition
\[
  {\mathbf f}_\alpha(x,y) = {\mathbf g}_\alpha(x,y) + t_\alpha \qquad (x,y \in G, \, \alpha \in {\mathbb A}),
\]
where $t_\alpha \geq 0$, and  ${\mathbf g}_\alpha |_G$ is a positive definite function associated with a representation weakly contained in $\hat{G} \setminus \{ 1_G \}$. As in the proof of Propositon \ref{oleg3} ---
possibly passing to a subnet ---, we conclude that $t_\alpha \to 1$ and consequently that $\mathbf{g}_\alpha(e,e) = \| \mathbf{g}_\alpha \| \to 0$. It follows that ${\mathbf f}_\alpha \to 1_{G \times G}$ in the norm of $B(G \times G)$, which is impossible by
Definition \ref{appind}(b).
\end{proof}
\begin{remarks}
\item Even though Theorem \ref{oleg4} makes $\SL(3,\comps)$ the prime suspect for a non-biflat, locally compact group, we would like to emphasize that the theorem does not prove this: All it shows is that Proposition \ref{biflatprop}
cannot be used to establish its biflatness. It is possible that our definition of an approximate indicator is too restrictive. Maybe one should consider nets, not in the Fourier--Stieltjes algebra, but consisting of completely bounded
multipliers of the Fourier algebra. Such multipliers are studied in \cite{Spr2}.
\item We had remarked earlier that $\SL(3,\comps)$ cannot be continuously embedded into any amenable, locally compact group. Combining Theorem \ref{biflatthm} and Theorem \ref{oleg4}, we see that $\SL(3,\comps)$ cannot even be
continuously embedded into a $[\operatorname{QSIN}]$-group.
\end{remarks}
\renewcommand{\baselinestretch}{1.0}
\dated
\vfill
\begin{tabbing}
{\it Second author's address\/}: \= Department of Mathematical and Statistical Sciences \kill
{\it First author's address\/}: \> Chair of Higher Mathematics \\
\> Obninsk Institute of Nuclear Power Engineering \\
\> Studgorodok-1 \\
\> 249040 Obninsk \\
\> Russia \\[\medskipamount]
{\it E-mail\/}: \> {\tt aoleg@krona.obninsk.ru} \\[\bigskipamount]
{\it Second author's address\/}: \> Department of Mathematical and Statistical Sciences \\
\> University of Alberta \\
\> Edmonton, Alberta \\
\> Canada T6G 2G1 \\[\medskipamount]
{\it E-mail\/}: \> {\tt vrunde@ualberta.ca} \\[\medskipamount]
{\it URL\/}: \> {\tt http://www.math.ualberta.ca/$^\sim$runde/} \\[\bigskipamount]
{\it Third author's address\/}: \> Department of Mathematics \\
\> Texas A \& M University \\
\> College Station, TX 77843-3368 \\
\> USA \\[\medskipamount]
{\it E-mail\/}: \> {\tt spronk@math.tamu.edu} \\[\medskipamount]
{\it URL\/}: \> {\tt http://www.math.tamu.edu/$^\sim$spronk/}
\end{tabbing}
\end{document}

%% file: format.tex
\renewcommand{\baselinestretch}{1.2}
\newcommand{\dated}{\mbox{} \hfill {\small [{\tt \today}]}}

%% file: mathdefs.tex
\usepackage{amsmath,amssymb,amscd}
\newcommand{\pf}[1]{\trivlist \item[\hskip\labelsep\it #1\ ]}
\newcommand{\varpf}[1]{\trivlist \item[\hskip\labelsep\sc #1:]}
\newcommand{\qedbox}{$\rlap{$\sqcap$}\sqcup$}
\newcommand{\qed}{\qquad \qedbox \endtrivlist}
\newcommand{\varqed}{\hfill \rule{0.6em}{0.6em} \endtrivlist}
\newenvironment{proof}{\pf{Proof}}{\qed}

\newenvironment{remark}{\pf{Remark}}{\endtrivlist}
\newenvironment{remarks}{\pf{Remarks} 
   \begin{enumerate}}{\end{enumerate} \endtrivlist}
\newenvironment{example}{\pf{Example}}{\endtrivlist}
\newenvironment{examples}{\pf{Examples} 
   \begin{enumerate}}{\end{enumerate} \endtrivlist}
\newenvironment{items}{
  \begin{enumerate} 
                    
  }{\end{enumerate}}
\newenvironment{alphitems}{
  \begin{enumerate} 
                    
  }{\end{enumerate}}
\newenvironment{keywords}{\noindent\small {\it Keywords\/}:}{\vskip 4pt}
\newenvironment{classification}{\noindent\small 2000 {\it Mathematics Subject
Classification\/}:}{\vskip 12pt}

\newcommand{\defiff}{\quad:\Longleftrightarrow\quad}
\newcommand{\comps}{{\mathbb C}}
\newcommand{\reals}{{\mathbb R}}

\newcommand{\free}{{\mathbb F}}

\newcommand{\void}{\varnothing}
\newcommand{\tensor}{\otimes}
\newcommand{\ttensor}{\tilde{\otimes}}
\newcommand{\Tensor}{\hat{\otimes}}

\newcommand{\cstar}{{C^\ast}}
\newcommand{\wstar}{{W^\ast}}

\newcommand{\id}{{\mathrm{id}}}

\newcommand{\A}{{\mathfrak A}}

\newcommand{\M}{{\mathfrak M}}

\newcommand{\VN}{\operatorname{VN}}
\newcommand{\SL}{\operatorname{SL}}
\newcommand{\SIN}{\operatorname{SIN}}

\newcommand{\supp}{{\operatorname{supp}}}

\newcommand{\varcl}[1]{\overline{#1}}

%% file: thmdefs.tex
\newtheorem{theorem}{Theorem}[section]
\newtheorem{lemma}[theorem]{Lemma}
\newtheorem{corollary}[theorem]{Corollary}
\newtheorem{proposition}[theorem]{Proposition}
\newtheorem{df}[theorem]{Definition}
\newenvironment{definition}{\begin{df} \rm}{\end{df}}